
\baselineskip=14pt
\parskip=10pt
\def\halmos{\hbox{\vrule height0.15cm width0.01cm\vbox{\hrule height
  0.01cm width0.2cm \vskip0.15cm \hrule height 0.01cm width0.2cm}\vrule
  height0.15cm width 0.01cm}}

\magnification=\magstephalf

\def\1{{\overline{1}}}
\def\2{{\overline{2}}}
\parindent=0pt
\overfullrule=0in

\def\frac#1#2{{#1 \over #2}}

\bf
\centerline
{
Sketch of a Proof of an Intriguing Conjecture of Karola M\'esz\'aros and Alejandro Morales
}
\centerline
{
Regarding the Volume of the $D_n$ Analog of the Chan-Robbins-Yuen Polytope
}
\centerline
{
(Or: The Morris-Selberg Constant Term Identity Strikes Again!)
}
\rm
\bigskip
\centerline
{\it By  Doron ZEILBERGER}
\bigskip
\qquad \qquad \qquad \qquad 
{\it To Dick Askey (b. June 4, 1933): from a $4^3$-year-old to a $3^4$-year-old, and thanks for preaching the importance of Constant Term Identities!}

\bigskip

{\bf Disclaimer:} Some of the steps below (in particular, the ``change of variable'' in contour-integrals)
require `rigorous' justification, that I am sure could be easily supplied by a skilled analyst. 

Recall that for any rational function $f(z)$ of a variable $z$,
$CT_{z} f(z)$ is the coeff. of $z^0$ in the formal Laurent expansion of $f(z)$ (that always exists!).

Karola M\'esz\'aros and Alejandro Morales have recently made the following intriguing conjecture.

{\bf Conjecture} ([MeMora], Conj. 7.12, also presented by Morales[Mora] at Stanley@70)
$$
CT_{x_n} \, CT_{x_{n-1}} \, \dots \, CT_{x_1} \prod_{i=1}^n x_i^{-1} (1-x_i)^{-2} \, \prod_{1 \leq i < j \leq n} \,  (x_j - x_i)^{-1} (1-x_j - x_i)^{-1}
\quad = \quad
2^{n^2} \prod_{k=1}^{n} Cat(k) \quad ,
$$
where $Cat(k)$ are the ubiquitous  Catalan numbers $(2k)!/(k!(k+1)!)$ (that Igor Pak believes are better than primes for searching for ETI!)

The present conjecture is a $D_n$-analog of a conjecture made in [CRY], and proved in [Z], using the
versatile Morris-Selberg Constant Term Identity ([Morr], restated in [Z]):
$$
CT_{x_n} \dots CT_{x_1} \,\,
\prod_{i=1}^{n} (1-x_i)^{-a} \prod_{i=1}^{n} x_i^{-b} 
\prod_{1 \leq i < j \leq n}
(x_j - x_i)^{-2c}
=S_n(a,b,c) \quad,
\eqno(Chip)
$$ 
where
$$
S_n(a,b,c)\, := \, {{1} \over {n!}} \prod_{j=0}^{n-1}
{{\Gamma(a+b+(n-1+j)c)\Gamma(c)} \over
{\Gamma(a+jc)\Gamma(c+jc)\Gamma(b+jc+1)}} \quad .
$$
By using Cauchy's theorem, this is equivalent to
$$
(\frac{1}{2\pi i})^n \int_C \,\,
\prod_{i=1}^{n} (1-x_i)^{-a} \prod_{i=1}^{n} x_i^{-b-1} 
\prod_{1 \leq i < j \leq n}
(x_j - x_i)^{-2c} \prod_{i=1}^n d\,x_i \, = \,S_n(a,b,c) \quad ,
\eqno(Atle)
$$ 
where {\it now} $a,b,c$ can be {\it any} real numbers (with obvious conditions to ensure convergence) and $C$ is {\it any}
multi-contour in $n$-dimensional complex space, that is far enough from the origin.

The M\'esz\'aros-Morales conjecture is the special case $a=2$, $c=\frac{1}{2}$ of the following fact.

\vfill\eject

{\bf Fact} : Let $a$ be a positive integer, and $c$ a positive half-integer, then
$$
CT_{x_n} \, CT_{x_{n-1}} \, \dots \, CT_{x_1}  \,\,
\prod_{i=1}^n x_i^{-(a-1)} (1-x_i)^{-a} \, \prod_{1 \leq i < j \leq n} \, (x_j - x_i)^{-2c} (1-x_j - x_i)^{-2c}
$$
$$
= \, 2^{2cn(n-1)+2(a-1)n} \, \cdot \, S_n(a \, , \, -\frac{1}{2} \, , \, c) \quad .
$$

{\bf Proof:} Converting the iterated constant-terms to a multi-contour-integral, we have to evaluate
$$
(\frac{1}{2\pi i})^n \int_C \,
\prod_{i=1}^n x_i^{-a} (1-x_i)^{-a} \, \prod_{1 \leq i < j \leq n} \,  (x_j - x_i)^{-2c} (1-x_j - x_i)^{-2c}
 \prod_{i=1}^n d\,x_i \quad .
$$
Now make the {\it change of variables}:
$$
x_i= \frac{1-z_i}{2} \quad, \quad ( 1 \leq i \leq n) \quad ,
$$
getting that our multi-integral equals
$$
(\frac{1}{2\pi i})^n \, \cdot \, 2^{2an} \cdot 2^{ 2cn(n-1)} (-\frac{1}{2})^n  \int_{C'} \,
\prod_{i=1}^n (1-z_i)^{-a} (1+z_i)^{-a} \, \prod_{1 \leq i < j \leq n} \, (z_i - z_j)^{-2c} (z_i+z_j)^{-2c}
 \prod_{i=1}^n d\,z_i
$$
$$
= \, 
(\frac{1}{2\pi i})^n \, \cdot \, (-1)^n 2^{ 2an+ 2cn(n-1) -n} \int_{C'} \,\,
  \prod_{i=1}^n (1-z_i^2)^{-a} \, \prod_{1 \leq i < j \leq n} \, (z_i^2 - z_j^2)^{-2c} \prod_{i=1}^n d\,z_i \quad ,
$$
where $C'$ is some other multi-contour.
Now is time for {\it yet another} {\bf change of variable}
$$
w_i=z_i^2 \quad, \quad ( 1 \leq i \leq n) \quad ,
$$
and we have
$$
\prod_{i=1}^n d\,z_i = (\frac{1}{2})^n \prod_{i=1}^n \frac{d\,w_i}{w_i^{1/2}} \quad,
$$
getting that our multi-integral is
$$
(-1)^{n+n(n-1)c} (\frac{1}{2\pi i})^n \, 2^{ 2an+ 2cn(n-1) -2n} \int_{C''} \,
  \prod_{i=1}^n (1-w_i)^{-a}  \prod_{i=1}^n w_i^{-1/2} \, \prod_{1 \leq i < j \leq n} \, (w_j - w_i)^{-2c} \prod_{i=1}^n d\,w_i \quad ,
$$
for yet another multi-contour, $C''$, and thanks to $(Atle)$ this is
$$
 2^{2cn(n-1)+2(a-1)n} \, \cdot \, S_n(a \, , \, -\frac{1}{2} \, ,\, c) \quad ,
$$
times $(-1)^{n+n(n-1)c}$, and I am sure that this annoying extra sign could be explained away. \halmos

\vfill\eject

{\bf References}

[CRY] Clara S. Chan, David P. Robbins, and David S. Yuen, {\it On the volume of a certain polytope},
Experimental Mathematics {\bf 9}(2000), 91-99. \quad {\tt http://arxiv.org/abs/math/9810154} \quad .

[MeMora] Karola M\'esz\'aros and Alejandro H. Morales, {\it Flow polytopes of signed graphs and the Kostant partition function},
{\tt http://arxiv.org/abs/1208.0140}, to appear in  International Mathematics Research Notices.

[Mora]  Alejandro H. Morales, {\it Open Problems session}, Stanley@70, June 23, 2014, ca. 4:50-5:00pm,
{\tt http://math.mit.edu/stanley70/Site/Program.html} \quad . \hfill\break
Slides are available here: {\tt http://math.mit.edu/stanley70/Site/Slides/Morales.pdf}

[Morr] William (``Chip'') G. Morris, ``{\it Constant Term Identities for Finite and Affine Root Systems:
Conjectures and Theorems}'', PhD thesis, University of Wisconsin-Madison, 1982. [Advisor: Richard Askey].

[Z] Doron Zeilberger, {\it Proof of a Conjecture of Chan, Robbins, and Yuen},
Elec. Trans. of Numerical Analysis (ETNA), {\bf 9}(1999), 147-148.\hfill\break
{\tt http://www.math.rutgers.edu/\~{}zeilberg/mamarim/mamarimhtml/cry.html}.

\bigskip
\hrule
\bigskip
Doron Zeilberger, Department of Mathematics, Rutgers University (New Brunswick), Hill Center-Busch Campus, 110 Frelinghuysen
Rd., Piscataway, NJ 08854-8019, USA. \hfill\break
url: {\tt http://www.math.rutgers.edu/\~{}zeilberg/}   
\quad . \hfill\break
Email: {\tt zeilberg at math dot rutgers dot edu}   \quad .
\bigskip
\hrule
\bigskip
{\bf EXCLUSIVELY PUBLISHED IN THE PERSONAL JOURNAL OF SHALOSH B. EKHAD and DORON ZEILBERGER} 
{\tt http://www.math.rutgers.edu/\~{}zeilberg/pj.html} and {\tt arxiv.org} \quad .
\bigskip
\hrule
\bigskip
July 9, 2014
\end